\documentclass[11pt]{article}
\usepackage{amsmath,amssymb}
\usepackage{latexsym,bm}
\usepackage{indentfirst}
\usepackage{ulem}
\usepackage{CJK}

\begin{document}
\begin{CJK*}{GBK}{song}

\title{{\bf On the integrality of the elementary symmetric
functions of $1, 1/3, ..., 1/(2n-1)$}
\thanks{ Research was supported partially by National Science Foundation
of China Grant \#10971145 and by the Ph.D. Programs Foundation of Ministry
of Education of China Grant \#20100181110073}}
\author{Chunlin Wang \\
{\it Mathematical College, Sichuan University, Chengdu 610064, China}\\
{\small wdychl@126.com, chun.lin.wang@hotmail.com}\\
\\
Shaofang Hong
\thanks{Corresponding author} \\
{\it Yangtze Center of Mathematics, Sichuan University, Chengdu 610064, China}\\
{\small sfhong@scu.edu.cn, s-f.hong@tom.com, hongsf02@yahoo.com}}
\date{}
\maketitle


\begin{abstract}
Erd\H{o}s and Niven proved that for any
positive integers $m$ and $d$, there are only finitely many positive
integers $n$ for which one or more of the elementary symmetric
functions of $1/m,1/(m+d), ..., 1/(m+nd)$ are integers. Recently, Chen and Tang
proved that if $n\ge 4$, then none of the elementary symmetric
functions of $1,1/2, ..., 1/n$ is an integer. In this paper, we show
that if $n\ge 2$, then none of the elementary symmetric functions of $1, 1/3,
..., 1/(2n-1)$ is an integer.
\end{abstract}

\section*{1. Introduction}
A well-known result in number theory says that for any
positive integers $m, d$, if $n>1$, then the harmonic sum
$\sum_{i=1}^n\frac{1}{m+id}$ is not an integer.
In 1946, Erd\H{o}s and Niven [1] proved
that there are only finitely many integers $n$ for which one or more
of the elementary symmetric functions of $1, 1/2, ..., 1/n$ are
integers, and they mentioned that by a similar argument one could
acquire the same result for the elementary symmetric function of
$1/m, 1/(m+d), ..., 1/(m+nd)$ for any given positive integers $m$ and $d$.
Recently, Chen and Tang [3] proved
that if $n\geq4$, then none of the elementary symmetric functions of
$1, 1/2, ..., 1/n$ is not an integer. It is an interesting question to determine
all finite arithmetic progressions $\{m+di\}_{i=0}^{n}$ such that
one or more elementary symmetric functions of $1/m, 1/(m+d), ..., 1/(m+nd)$
are integers.

In this paper, we consider the finite arithmetic progression
$\{1+2i\}_{i=0}^{n-1}$. Throughout, we let $S_k(n)$ denote the $k$-th
elementary symmetric functions of $1, 1/3, ..., 1/(2n-1)$. That is,
$$
S_k(n):=\sum_{0\le i_1<...<i_k\le
n-1}\prod_{j=1}^k\frac{1}{(1+2i_j)}.
$$
We will show that for all integers $n>1\ {\rm and}\ 1\leq k\leq n$,
$S_k(n)$ is not an integer. See Theorem 3.1 below.
The paper is organized as follows. In
Section 2, we show several lemmas which are needed for the proof of
the main result. In the last section, we give the main result.

\section*{2. Several lemmas}

In the present section, we show some preliminary lemmas which are
needed for the proof of our main result. As usual, let $\pi (x)$
denote the number of primes no more than $x$.
We begin with a known result.\\
\\
{\bf Lemma 2.1.}  [2] {\it One has
$$\pi(x)<\frac{x}{\log x-1-(\log x)^{-1/2}} \ \text{for all}\ x\geq6$$
and
$$\pi(x)>\frac{x}{\log x-1+(\log x)^{-1/2}} \ \text{for all}\
x\geq59.$$}\\
\\
{\bf Lemma 2.2.} {\it For any integer $k\geq1$, we have
$$
S_{k}(k+1)=\frac{(k+1)^2}{\prod_{i=0}^k(1+2i)}, S_{k}(k+2)
=\frac{(k+1)(k+2)(3k^{2}+11k+9)}{6\prod_{i=0}^{k+1}(1+2i)}
$$
and
$$
S_{k}(k+3)=\frac{(k+1)(k+2)(k+3)^{2}(k^{2}+5k+5)}{6\prod_{i=0}^{k+2}(1+2i)}.
$$
}
\par\noindent {\it Proof.} Since
$$
\sum_{i=0}^k(1+2i)=(k+1)^{2},
\sum_{0\leq i<j\leq k+1}(1+2i)(1+2j)=\frac{1}{6}(k+1)(k+2)(3k^{2}+11k+9)
$$
and
$$
\sum_{0\leq i<j<l\le k+2}(1+2i)(1+2j)(1+2l)=\frac{1}{6}(k+1)(k+2)(k+3)^{2}(k^{2}+5k+5).
$$
the desired formulae follow immediately. So Lemma 2.2 is proved. \hfill$\Box$ \\
\\
{\bf Lemma 2.3.} {\it Let $k$ and $n$ be positive integers
such that
$$e\bigg(\frac{1}{2}\log(2n-1)+1\bigg)\le k\le n.$$
Then $S_{k}(n)$ is not an integer.}\\
\\
{\it Proof.} First by the multi-nomial expansion theorem, we
get
$$S_{k}(n) \leq \frac{1}{k!}\bigg(\sum\limits_{i=0}^{n-1}\frac{1}{1+2i}\bigg)^{k}.$$

On the one hand, one has
$$\sum\limits_{i=0}^{n-1}\frac{1}{1+2i}<1+\int_{0}^{n-1}
\frac{1}{1+2x}dx=\frac{1}{2}\log(2n-1)+1.$$
On the other hand, we have
$$\log k!=\sum_{i=2}^{k}\log i>\int_{1}^{k}\log x {\rm d}x
>k\log k-k>k\log(\frac{1}{2}\log(2n-1)+1).$$
So from the above inequalities, we deduce that
$$
\bigg(\sum\limits_{i=0}^{n-1}\frac{1}{1+2i}\bigg)^{k}<k!.
$$
In other words, $S_{k}(n)<1$ if $n \geq k \geq e(\frac{1}{2}\log(2n-1)+1)$.
This ends the proof of Lemma 2.3. \hfill$\Box$ \\
\\
{\bf Lemma 2.4.} {\it Let $k$ and $n$ be positive integers such that $1<k\le n$.
Suppose that there exists an odd prime $p>2k+6$ satisfying that
$$\frac{n}{k+3}<p\leq\frac{n}{k}$$
and
$$p\nmid(3k^{2}+11k+9)(k^{2}+5k+5).$$
Then $S_{k}(n)$ is not an integer.}\\
\\
{\it Proof.} First of all, we can easily check that the following
identity holds:
$$
S_{k}(n)=\sum_{0\leq i_1<\cdots<i_k\leq[\frac{n}{p}]+t}\prod_{j=1}^{k}\frac{1}{p(1+2i_j)}+
\sum_{\substack{0\leq i_1<\cdots<i_k\leq n-1 \\
\exists j\ {\rm s.t.\ }p\nmid(1+2i_j)}}\prod_{j=1}^{k}\frac{1}{(1+2i_j)}, \eqno(2.1)
$$
where $t=-1\ \text{if}\ p(1+2[\frac{n}{p}])>2n-1$, and $t=0$
otherwise.

Since $p>2k+6\ {\rm and}\ p>\frac{n}{k+3}$, we have $p>\sqrt{2n-1}$.
It infers that $v_p(1+2i)\le 1$ for $0\leq i\leq n-1$, where $v_p$
denotes the $p$-adic valuation on $\mathbb{Q}$. We then derive from (2.1) that
$$
S_k(n)=\frac{1}{p^k}S_k\bigg(\bigg[\frac{n}{p}\bigg]+t+1\bigg)+\frac{a}{p^{k-1}b} \eqno(2.2)
$$
for some positive integers $a$ and $b$ with $p\nmid b$. Note that
$k\leq[\frac{n}{p}]+t+1\leq k+3$. But $p>2k+6$ and
$p\nmid(3k^{2}+11k+9)(k^{2}+5k+5)$. Then by Lemma 2.2 we obtain that
$v_p(S_k([\frac{n}{p}]+t+1))=0$.

Now using (2.2), we can get that $v_{p}(S_{k}(n))=-k<0.$ Therefore
$S_{k}(n)$ is not an integer as desired. The proof of Lemma 2.4 is
complete. \hfill$\Box$

\section*{3. The main result}

In this section, we give the main result of this paper.\\
\\
{\bf Theorem 3.1.} {\it For any integers $n>1$ and $k$ with
$1\leq k\leq n$, $S_{k}(n)$ is not an integer.}\\
\\
{\it Proof.} When $k=1$, it is known that for any integer $n>1$,
$\sum_{i=0}^{n-1}\frac{1}{1+2i}$ is not an integer (see, for
example, [1]). So Theorem 3.1 is true when $k=1$. In what follows we
let $k\geq2$.

By Lemma 2.3, we know that $S_{k}(n)$ is not an integer if
$e(\frac{1}{2}\log(2n-1)+1)\le k\le n$. In the following we assume
that $2\le k<e(\frac{1}{2}\log(2n-1)+1)$.

First we let $n\ge 23000$. Claim that there is a prime number
$p>2k+6$ such that $\frac{n}{k+3}<p\leq\frac{n}{k}$ and
$p\nmid(3k^{2}+11k+9)(k^{2}+5k+5)$. It then follows immediately from
the claim and Lemma 2.4 that $S_{k}(n)$ is not an integer for all
$2\le k<e(\frac{1}{2}\log(2n-1)+1)$ if $n\ge 23000$. It remains to
show the claim which will be done in the following.

First we prove $\pi(\frac{n}{k})>\pi(\frac{n}{k+3})$. By Lemma 2.1, it
suffices to show that
$$
\frac{n/k}{\log(n/k)-1+(\log(n/k))^{-1/2}}>\frac{n/(k+3)}
{\log(n/(k+3))-1-(\log(n/(k+3)))^{-1/2}}.
$$
This is equivalent to
$$
k\log\bigg(1+\frac{3}{k}\bigg)+3+k\bigg(\log\frac{n}{k}\bigg)^{-1/2}+(k+3)
\bigg(\log\frac{n}{k+3}\bigg)^{-1/2}<3\log\frac{n}{k+3}. \eqno(3.1)
$$

Since $k\log(1+\frac{3}{k})<3$ and
$\log\frac{n}{k}>\log\frac{n}{k+3}$, in order to show that (3.1)
holds, it is enough to prove that the following inequality is true:
$$
6+(2k+3)\bigg(\log\frac{n}{k+3}\bigg)^{-1/2}<3\log\frac{n}{k+3}. \eqno(3.2)
$$

Define a real function $f(x)$ by
$$f(x):=x^{0.3}-\frac{e}{2}\log(2x-1)-e-3.$$
Then one can easily check that $f(23000)>0$ and
$$xf'(x)=0.3x^{0.3}-\frac{ex}{2x-1}>0$$ for all $x\ge 23000$. We can derive that
$f(x)>0$ for all $x\ge 23000$. But $k<e(\frac{1}{2}\log(2n-1)+1)$.
So for $n\ge 23000$, we have
$$\frac{n}{k+3}>\frac{n}{\frac{e}{2}\log(2x-1)+e+3}>n^{0.7}.$$
Thus to prove (3.2), it is sufficient to show the following
inequality
$$6+(e\log(2n-1)+2e+3)(\log n^{0.7})^{-1/2}\leq3\log n^{0.7}, n\ge 23000,$$
which is equivalent to
$$6\times0.7^{1/2}(\log n)^{1/2}+e\log(2n-1)+2e+3\le 3
\times0.7^{3/2}(\log n)^{3/2}, n\ge 23000. \eqno(3.3)
$$
Let $t=(\log n)^{1/2}$. Then one find that for $t\ge 3.10$, we have
$$
6\times0.7^{1/2}t+et^2+e\log2+2e+3\le 3 \times0.7^{3/2}t^3,
$$
from which (3.3) follows immediately. Hence (3.2) is proved and so we have
$\pi(\frac{n}{k})>\pi(\frac{n}{k+3})$ for
$k<e(\frac{1}{2}\log(2n-1)+1)$ if $n\geq23000$.

Consequently, we prove that $\frac{n}{k+3}>3k^{2}+11k+9$ for
$k<e(\frac{1}{2}\log(2n-1)+1)$ if $n\geq23000$. Evidently we have
$n>\frac{1}{2}e^{(2k/e)-2}$ since $k<e(\frac{1}{2}\log(2n-1)+1)$. It
is easy to show that
$$\frac{1}{2}e^{(2x/e)-2}>g(x):=(x+3)(3x^{2}+11x+9)$$
for all $x\ge 17.3$. Let $h(x):=e(\frac{1}{2}\log(2x-1)+1)$. Then
$h(n)\ge 17.3$ if $n\geq23000$ and $h(n)>k$. It follows that
$$
n=\frac{1}{2}e^{\frac{2h(n)}{e}-2}+\frac{1}{2}
>\frac{1}{2}e^{\frac{2h(n)}{e}-2}>g(h(n))>g(k).
$$
Namely, we have $\frac{n}{k+3}>3k^{2}+11k+9$ for
$k<e(\frac{1}{2}\log(2n-1)+1)$ if $n\ge 23000$.

Since $\pi(\frac{n}{k})>\pi(\frac{n}{k+3})$, there is a prime number
$p$ satisfying $\frac{n}{k+3}<p\leq\frac{n}{k}$. But
$\frac{n}{k+3}>3k^{2}+11k+9$. Thus $p>2k+6$ and
$p\nmid(3k^{2}+11k+9)(k^{2}+5k+5)$. Hence the claim is proved.

Now we treat the remaining case: $n<23000$. Since
$k<e(\frac{1}{2}\log(2n-1)+1)$, we have $k<18$ and
$n>\frac{1}{2}e^{2k/e-2}+\frac{1}{2}$.

If $12\leq k\leq17$, then $\frac{1}{2}e^{2k/e-2}>2(k+3)^2$. This
implies that $\frac{n}{k+3}>2k+6$. We can check by computer that for
every integer $n\in(\frac{1}{2}e^{2k/e-2},23000)$, there is a prime
number $p$ such that
$$\frac{n}{k+3}<p\leq\frac{n}{k} \ {\rm and}\ p\nmid(k^2+5k+5)(3k^2+11k+9).$$
Hence by Lemma 2.4, we know that $S_k(n)$ is not an integer for
$12\leq k\leq17$ and $n<23000$.

If $2\leq k\leq11$, then $400>2(k+3)^2$. It can be checked by
computer that for every integer $400\leq n<23000$, there is a prime number
$p$ such that
$$
\frac{n}{k+3}<p\leq\frac{n}{k}\ {\rm and} \ p\nmid(k^2+5k+5)(3k^2+11k+9).
$$
Note that for the above prime $p$,
we have $p>\frac{n}{k+3}>2k+6$. Then Lemma 2.4 tells us that $S_k(n)$
is not an integer if $2\leq k\leq11$ and $n\ge 400$.

For the remaining case $n<400$, we can verify by using Maple 12 that
$S_k(n)$ is not an integer for any integers $2\le k\le 11$ and $n<400$.

This completes the proof of Theorem 3.1. \hfill$\Box$

\section*{Reference}
\par\noindent [1] P. Erd\H{o}s and I. Niven, Some properties of partial sums of the
 harmonic series, {\it Bull. Amer. Math Soc.} 52 (1946), 248-251.\\
\par\noindent [2] L. Panaitopol, Inequalities concerning the function $\pi(x)$: Applications,
{\it Acta Arith.} 94 (2000), 373-381.\\
\par\noindent [3] Y. Chen and M. Tang, On the elementary symmetric functions of
$1, 1/2, ...,$ $1/n$, arXiv:1109.1142.


\end{CJK*}
\end{document}